\documentclass{article}

\usepackage{amssymb}
\usepackage{latexsym}
\usepackage{amsthm}
\usepackage{enumerate}
\usepackage{epsfig}
\usepackage{color}
\usepackage{float}
\usepackage{subfigure}
\usepackage{amsmath}
\usepackage{ifthen}
\usepackage{makeidx}
\usepackage{lscape}
\usepackage{setspace}
\usepackage{amscd}
\usepackage{xspace}
\usepackage[nobysame,alphabetic]{amsrefs}
\usepackage{srcltx}
\usepackage{fancyvrb}

\newboolean{draft}
\setboolean{draft}{false}

\newcommand{\mylabel}[1]{\label{#1}\ifthenelse{\boolean{draft}}{{\bf [ #1 ]}}{}}
\newcommand{\pref}[1]{Proposition \ref{prop:#1}\xspace}

\makeatletter
\newcommand*{\mydef}[2]{%
    \def#1{#2}%
    \ifthenelse{\boolean{draft}}{{\noindent\ttfamily\string#1=\expandafter\strip@prefix\meaning#1\\}}{}%
}
\makeatother

\def\a{\alpha}

\def\C{\mathcal{C}}
\def\E{\mathbb{E}}
\def\g{\gamma}

\def\D{\Delta}

\def\l{\lambda}

\def\R{\mathbb{R}}

\def\S{\Sigma}

\def\Z{\mathbb{Z}}
\def\d{\partial}
\def\cross{\times}

\def\e{\epsilon}
\def\half{\frac{1}{2}}
\def\t{\tau}

\def\P{\mathbb{P}}
\def\le{\leqslant}
\def\ge{\geqslant}
\def\<{\langle}
\def\>{\rangle}

\def\PMF{\mathcal{PMF}}

\def\fmin{\mathcal{F}_{min}}
\def\Grel{\widehat G}

\def\Gbar{\overline G}

\def\rmu{\widetilde \mu}
\def\rnu{\widetilde \nu}

\def\pA{pseudo-Anosov\xspace}
\def\mcg{mapping class group\xspace}

\newcommand{\norm}[1]{\left| #1 \right|}

\newcommand{\dhat}[1]{d (#1)}
\newcommand{\nhat}[1]{\dhat{1, #1}}

\newcommand{\gp}[3]{(#2\cdot#3)_{#1}}
\newcommand{\mun}[1]{\mu_{#1}}
\newcommand{\rmun}[1]{\widetilde \mu_{#1}}

\newcommand\pic[3]{
\begin{figure}[H] \begin{center} 
\epsfig{file=#1, height=#2pt} 
\end{center} 
\caption{#3} 
\label{pic:#1}
\end{figure}
}

\newtheorem{theorem}{Theorem}[section]
\newtheorem{lemma}[theorem]{Lemma}

\newtheorem{proposition}[theorem]{Proposition}

\theoremstyle{definition}
\newtheorem{properties}[theorem]{Properties}

\restylefloat{figure}

\begin{document}

\title{Exponential decay in the \mcg}
\author{Joseph Maher\footnote{email: joseph.maher@csi.cuny.edu}}
\date{\today}

\maketitle

\begin{abstract}
  We show that the probability that a non-pseudo-Anosov element arises
  from a finitely supported random walk on a non-elementary subgroup
  of the mapping class group decays exponentially in the length of the
  random walk. More generally, we show that if $R$ is a set of mapping
  class group elements with an upper bound on their translation
  lengths on the complex of curves, then the probability that a random
  walk lies in $R$ decays exponentially in the length of the random
  walk.
\end{abstract}

\tableofcontents

\section{Introduction}

Let $\S$ be a compact oriented surface of finite type, and let $G$ be
the mapping class group of $\S$, i.e. the group of orientation
preserving homeomorphisms of $\S$, up to isotopy.  Let $\mu$ be a
probability distribution on $G$ with finite support. A random walk on
$G$ is a Markov chain with transition probabilities $p(x, y) =
\mu(x^{-1}y)$.  We will always assume we start at the identity at time
zero, and we will write $w_n$ for the location of the random walk at
time $n$.  The probability distribution $\mu$ need not be symmetric,
but we shall always assume that the semi-group generated by the
support of $\mu$ contains a non-elementary subgroup of the mapping class
group.  A subgroup of the mapping class group is \emph{non-elementary}
if it contains a pair of \pA elements with distinct fixed points in
$\PMF$, Thurston's boundary for the mapping class group.  Rivin
\cite{rivin, rivin2} and Kowalski \cite{kowalski} showed that the
probability that a random walk on the mapping class group gives rise
to a pseudo-Anosov element tends to one, as long as the group
generated by the support of the mapping class group maps onto a
sufficiently large subgroup of $Sp(2g, \Z)$. Furthermore, they showed
that the probability that an element is \emph{not} pseudo-Anosov
decays exponentially in the length of the random walk. Malestein and
Souto \cite{ms} and Lubotzky and Meiri \cite{lm} extended this to the
Torelli subgroup, by considering the action of the Torelli group on
the homology of double covers of the surface.

In \cite{maher1} it was shown that the probability that a random walk
gives a \pA element tends to one for all non-elementary subgroups of
the mapping class, by considering the action of the mapping class
group on the complex of curves, but no information was obtained about
the rate of convergence. In this paper we show that the rate of
convergence is exponential; in fact, we show that for any constant
$B$, the probability that a random walk gives an element of
translation length at most $B$ on the complex of curves decays
exponentially in the length of the random walk; the rate of decay
depends on $B$. Furthermore, the argument given here avoids various
complications involving centralizers that arose in the earlier
approach.

We say the surface $\S$ is \emph{sporadic} if $\S$ is a sphere with at
most four punctures, or a torus with at most one puncture. The complex
of curves $\C(\S)$ is a simplicial complex, whose vertices consist of
isotopy classes of simple closed curves, and whose simplices are
spanned by disjoint simple closed curves. The mapping class group $G$
acts by simplicial isometries on the complex of curves, and Masur and
Minsky \cite{mm1} showed that an element is \pA if and only if its
translation length on $\C(\S)$ is positive. We will use Landau's ``big
$O$'' notation, so $O(g(x))$ denotes some function $f(x)$ such that
$f(x) \le C \norm{g(x)}$ for some constant $C > 0$, and for all $x$
sufficiently large.

\begin{theorem} \mylabel{theorem:main}
Let $G$ be the mapping class group of a non-sporadic surface of finite
type, and let $w_n$ be a random walk of length $n$ on $G$ generated by
a finitely supported probability distribution $\mu$, whose support
generates a non-elementary subgroup of the mapping class group. Then
for any constant $B > 0$, there is a constant $c < 1$ such that
\[ \P(\tau(w_n) \le B) \le O( c^n), \]
where $\tau(w_n)$ is the translation length of $w_n$ acting on the
complex of curves.
\end{theorem}

If the surface is sporadic, then the mapping class group is either
finite, or commensurable to $SL(2, \Z)$, and in the latter case the
result follows from the work of Rivin \cite{rivin, rivin2} or Kowalski
\cite{kowalski} on random walks on matrix groups.  Theorem
\ref{theorem:main} does not apply to the Torelli group of the genus
two surface, as this group is not finitely generated, as shown by
McCullough and Miller \cite{mcmi}. However, the results of
\cite{maher1} hold in this case, but with no rate of convergence
information.

There are two main steps, both of which use the improper metric on $G$
arising from its action on the complex of curves, which we shall
denote by $d(g, h) = d_{\C(\S)}(g x_0, h x_0)$, where $x_0$ is a
basepoint in the complex of curves $\C(\S)$; this is also known as a
relative metric on $G$. The first step is to show that the random walk
has a linear rate of escape in the relative metric, with exponential
decay for the proportion of sample paths making progress at lower
rate. The second is to consider the distribution of elements of
bounded translation length on the complex of curves. If $g$ is an
element of bounded translation length on the complex of curves, then
$g$ is conjugate to an element $s$, of bounded relative
length. Furthermore, if $v$ is chosen to be a shortest conjugating
element, then the path $vsv^{-1}$ is quasigeodesic, with quasigeodesic
constants depending only on $G$ and the bound on translation
length. This means that if a random walk $w_n$ is conjugate to an
element of bounded translation length, then if the first half of a
geodesic from $1$ to $w_n$ fellow travels with some geodesic from $1$
to $v$, then the second half of the geodesic from $1$ to $w_n$ fellow
travels with a translate of a geodesic from $1$ to $v^{-1}$. This
fellow travelling condition is equivalent to the condition that the
pair $(w_n, w_n^{-1})$ lies in a certain neighbourhood of the diagonal
in $G \cross G$, and we show that the probability that this occurs
decays exponentially in the length of $w_n$.


We now give a brief summary of the organization of the paper. The
remainder of this section is devoted to a detailed outline of the
argument described in the previous paragraph. In Section
\ref{section:preliminaries} we introduce some standard definitions and
fix notation. In particular, we define subsets of $G$, called shadows,
and find upper bounds for the probability that a random walk lies in a
shadow.  In Section \ref{section:linear progress} we show the linear
progress result, with an exponential decay bound for the proportion of
paths making linear progress below some rate. Finally, in Section
\ref{section:translation length} we show that the fellow travelling
property described above is equivalent to the condition that the pair
$(w_n, w_n^{-1})$ lies in a certain neighbourhood of the diagonal in
$G \cross G$, consisting of unions of shadows, and we show that the
probability that a random walk lies in one of these neighbourhoods
decays exponentially in the length of $w_n$.

\subsection{Outline}

We will consider the action of the mapping class group on the complex
of curves, see Farb and Margalit \cite{fm} for an introduction to the
mapping class group.  The complex of curves $\C(\S)$ is a simplicial
complex, whose vertices are isotopy classes of simple closed curves,
and whose simplices are spanned by disjoint simple closed curves. The
complex of curves is finite dimensional, but not locally finite. We
will only need to consider distances between vertices in the curve
complex, and so we consider the $1$-skeleton of the complex of curves
to be a metric space $(\C(\S), d_{\C(\S)})$, by assigning every edge
to have length $1$. By abuse of notation, we will refer to this as a
metric on the curve complex.  The mapping class group acts on the
complex of curves by simplicial isometries, and a choice of basepoint
$x_0$ in the complex of curves determines a map from $G$ to $\C(\S)$
defined by $g \mapsto g(x_0)$. We may therefore define an improper
metric on the mapping class group by
\[ d(g, h) = d_{\C(\S)}(g x_0, h x_0). \]
We emphasize that throughout this paper the metric $d$ will always refer
to this improper metric induced from the action of the mapping class
group on the complex of curves, and never a proper word metric on $G$
with respect to a finite generating set. However, the metric $d$ is
quasi-isometric to a word metric on $G$ with respect to an infinite
generating set, also known as a \emph{relative metric}, formed by
starting with a finite generating set and adding subgroups which
stabilize vertices in the complex of curves which correspond to
distinct orbits under the action of the mapping class group. Masur and
Minsky \cite{mm1} showed that the curve complex is Gromov hyperbolic,
and we will write $\d G$ for the Gromov boundary of $G$, and $\Gbar$
for $G \cup \d G$.

A random walk of length $n$ on $G$ is a product of $n$ independent
identically $\mu$-distributed random variables $s_i$, which we shall
call the \emph{steps} of the random walk, so $w_n = s_1 s_2 \ldots
s_n$, and $w_n$ is distributed as $\mun{n}$, the $n$-fold convolution
of $\mu$ with itself.  A random walk converges to the Gromov boundary
almost surely, so this gives a hitting measure, known as
\emph{harmonic measure} on $\d G$, and which we shall denote by $\nu$.
We will need to estimate the probability that a random walk lies in
particular subsets of $\Gbar$. We now define a family of subsets of
$\Gbar$ which we shall call \emph{shadows}.  Recall that the Gromov
product of $x$ and $y$ with respect to $1$ is equal to the distance
from $1$ to a geodesic from $x$ to $y$, up to a bounded error which
only depends on $\delta$.  Given a real number $r$, we can use the
Gromov product to define the shadow of a point $x$ in $G$, which we
shall denote by $S_1(x, r)$,
\[ S_1(x, r) = \{y \in \Gbar \mid \gp{1}{x}{y} \ge r \}. \]
\pic{}{80}{A shadow of a point.}
We warn the reader that we use a different parameterization of shadows
than that used by other authors, for example, Blach\`ere,
Ha\"issinsky, and Mathieu, \cite{bhm} define their shadows $\mho_1(x,
r)$ to be $S_1(x, d(1, x) - r) \cap \d G$ in our notation. We show
that both the harmonic measures of shadows, and the $\mun{n}$-measures
of shadows, decay exponentially in $r$, i.e. there are constants $K$
and $c < 1$ such that $\nu(S_1(x, r)) \le c^r$ and $\mun{n}(S_1(x,
r)) \le K c^r$, for all $x$, $r$ and $n$.

In \cite{maher2}, we showed that a random walk makes linear progress
in the relative metric, almost surely, i.e. there is a constant $L >
0$ such that
\[ \P(d(1, w_n) \le Ln) \to 0 \text{ as } n \to \infty. \]
We need a stronger version of this result, which gives an exponential
decay bound for the rate of convergence. To be precise, we show:

\begin{theorem} \mylabel{theorem:exponential linear progress}
Let $G$ be the mapping class group of a non-sporadic surface of finite
type, and let $w_n$ be a random walk of length $n$ on $G$, generated
by a finitely supported probability distribution $\mu$, whose 
support generates a non-elementary subgroup of the
mapping class group. Then there are constants $L > 0$ and $c < 1$,
such that 
\[ \P \left( \dhat{1, w_n} \le Ln \right) \le O( c^n ), \]
where $d$ is the non-proper metric on the mapping class group arising
from its action on the complex of curves.
\end{theorem}

We now give a brief overview of the proof of Theorem
\ref{theorem:exponential linear progress}.  Consider taking the random
walk $k$ steps at a time, i.e. consider $w_{nk}$ instead of $w_n$,
which we shall refer to as the $k$-iterated random walk. We shall
write $w^k_n$ for $w_{nk}$, and the steps of the $k$-iterated random
walk are given by $s^k_n = s_{kn-k} s_{kn -k + 1} \ldots s_{kn}$.  The
increments of the walk, $s^k_n$, are all independent and identically
distributed, with distribution $\mun{k}$.  However, the distance
travelled at time $nk$, given by $d(1, w^k_{n})$ is not the sum of the
distances travelled at each step of the $k$-iterated walk, $d(w^k_{n},
w^k_{n+1})$, as there may be some ``backtracking'', illustrated
schematically in Figure \ref{pic:backtracks}.

\pic{}{100}{Steps of the iterated random walk.}

The distance $d(1, w^k_n)$ is the sum of the $d(w^k_{i-1}, w^k_{i})$
for $i \le n$, minus the total amount of backtracking. A key fact is
that the distribution of the amount of backtracking at time $i$ is
bounded above by an exponential function, and furthermore, the same
upper bound holds for all $i$, independent of the locations of the
random walk, or the amount of backtracking at other times, and also
independent of $k$, the number of steps for each segment of the
$k$-iterated random walk. We now explain why this is the case.  The
amount of backtracking can be estimated as follows. The size of the
backtrack from $w^k_{i-1}$ to $w^k_{i}$ is roughly the distance from
$w^k_{i-1}$ to the geodesic from $1$ to $w^k_{i}$. After applying the
isometry $w^k_{i-1}$, this is the same as the distance from $1$ to the
geodesic from $(w^k_{i-1})^{-1}$ to $(w^k_{i-1})^{-1} w^k_{i}$, as
illustrated in Figure \ref{pic:increment}.

\pic{}{70}{A single backtrack, after the isometry $w^k_{i-1}$.}

The point $(w^k_{i-1})^{-1} w^k_{i}$ is equal to $s^k_{i}$, and the
pair of random variables $( (w^k_{i-1})^{-1} , s^k_{i} )$ are
independent, and distributed as $\rmun{k(i-1)} \cross \mun{k}$, where
$\rmun{n}$ is the $n$-fold convolution of the reflected distribution
$\rmu(g) = \mu(g^{-1})$. If the backtrack is of length at least $r$,
then distance from $1$ to a geodesic from $(w^k_{i-1})^{-1}$ to
$s^k_{i}$ is at least $r$, up to bounded error, so in turn the Gromov
product $\gp{1}{(w^k_{i-1})^{-1}}{s^k_{i}}$ is at least $r$ up to
bounded error. Therefore $s^k_i$ lies in $S_1((w^k_{i-1})^{-1}, r -
K)$, for some $K$ which only depends on $\delta$. We show that both
the harmonic measure $\nu$, and the convolution measures $\mu_n$, of
shadows $S_1(x, r)$ are bounded above by a function which decays
exponentially in $r$, and furthermore, the upper bound function is
independent of both $x$ and $n$. Therefore, the probability that there
is a backtrack of size $r$ decays exponentially in $r$, independently
of $k$, and also independently of the locations of the random walk at
other times. In particular, the expected size of a backtrack is
bounded independently of $k$, so by choosing $k$ sufficiently large,
we can ensure that the expected value of each $k$-iterated step
$d(w^k_{i-1}, w^k_{i})$ is larger than the expected value of a
backtrack. Furthermore, applying standard Bernstein or
Chernoff-Hoeffding estimates for concentration of measures, we obtain
bounds for the probability that the sums of the first $n$ backtracks
and $k$-iterated steps deviate from their expected values, and these
bounds decay exponentially in $n$. This implies that the distance away
from the origin grows linearly at some rate, with exponential decay
for the proportion of paths making progress below this rate.

We now wish to show that the probability that $w_n$ is \pA tends to
$1$ exponentially quickly. The translation length of a group element
on the complex of curves is
\[ \tau(g) =  \lim_{n \to \infty} \tfrac{1}{n} d(1, g^n).  \]
Masur and Minsky \cite{mm1} showed that the \pA elements are precisely
those elements with non-zero translation length.  The translation
length of an element $g$ acting on the complex of curves is also
coarsely equivalent to the shortest length of any conjugate of $g$,
measured in the relative metric on the mapping class group.  In
\cite{maher1} we showed that the \mcg has \emph{relative conjugacy
  bounds}, i.e. there is a constant $K$ such that if two group
elements $a$ and $b$ are conjugate, then $a = vbv^{-1}$ for some
element $v$ whose length is bounded in terms of the lengths of $a$ and
$b$,
\[ d(1, v) \le K (d(1, a) + d(1, b) ). \]
We emphasize that the distance $d$ here is the relative or curve
complex distance on the mapping class group. Every group element $g$
corresponds to a point in $G$, but we may also think of $g$ as
representing some choice of geodesic from $1$ to $g$. We may therefore
think of a product of group elements as representing a path in $G$,
composed of concatenating various translates of geodesics representing
each element in the product. In particular, the word $vbv^{-1}$
corresponds to a path consisting of three geodesic segments.  As the
curve complex is $\delta$-hyperbolic, one may show that if $v$ is
chosen to be a conjugating word of shortest relative length, then the
path $vbv^{-1}$ is a quasigeodesic path in the curve complex, where
the quasigeodesic constants depend on the relative conjugacy bound
constant $K$, and the length of $b$.

If we choose $R$ to be a collection of group elements of conjugacy
length at most $B$, then every element $g \in R$ is equal to
$vsv^{-1}$, where $d(1, s) \le B$, and the paths $vsv^{-1}$ are
uniformly quasigeodesic over all elements of $R$. This implies that
the first half of the geodesic from $1$ to $g$ fellow travels with a
translate of the inverse of the second half of the geodesic from $1$
to $w_n$. In order to find an upper bound on the probability that this
occurs, it is convenient to express this fellow travelling property in
terms of the location of the pair $(g, g^{-1})$ in $G \cross G$. The
fact that $vsv^{-1}$ is quasigeodesic implies that $g \in S_1(v, r)$
and $g^{-1} \in S_1(v, r)$, where $r$ is equal to $\tfrac{1}{2} d(1,
g)$, up to an additive error which only depends on $\delta$ and the
quasigeodesic constants. We may extend the definitions of shadows to
subsets $U \subset G$ by setting $S_1(U, r)$ to be the union of all
$S_1(g, r)$ over all points $g \in U$. We may then extend the
definition of shadows to subsets $U \subset \Gbar \cross \Gbar$, by
setting $S_1(U, r)$ to be the union of all $S_1(g_1, r) \cross
S_1(g_2, r)$, over all $(g_1, g_2) \in U$. In particular, if a random
walk $w_n$ lies in $R$, then the pair $(w_n, w_n^{-1})$ lies in
$S_1(\D, r)$, a shadow of the diagonal $\D$ in $\Gbar \cross \Gbar$,
where $r$ is roughly $\tfrac{1}{2} d(1, w_n)$.  By the linear progress
result, we may assume that $r$ grows linearly in $n$, up to a set of
paths whose measure decays exponentially in $n$.  The distribution of
pairs $(w_n, w_n^{-1})$ is obviously not independent, as $w_n$
determines $w_n^{-1}$, but they are asymptotically independent, and
converge to $\nu \cross \rnu$.  In fact, we may approximate the
distribution of pairs $(w_{2n}, w_{2n}^{-1})$ by the distribution of
pairs $(w_n, w_{2n}^{-1}w_n)$. This is because as sample paths
converge to the boundary almost surely, it is probable that the the
point $w_n$ looks close to the point $w_{2n}$, as viewed from the
origin $1$, as illustrated in Figure \ref{pic:w2n} below.
\pic{}{80}{A path of length $2n$.}
Similarly, standing at $w_{2n}$ and looking back towards the origin
$1$, the point $1$ looks close to the midpoint of the path $w_n$. If
we apply the isometry $w_{2n}^{-1}$, this implies that $w_{2n}^{-1}$
and $w_{2n}^{-1} w_n$ look close together when viewed from $1$. We can
make this precise, and we show that the probability that $w_{2n}$ lies
in the shadow $S_1(w_n, d(1, w_n) - K)$ tends to one exponentially
quickly, for some $K$ which only depends on the constant of
hyperbolicity $\delta$. The same argument shows that the probability
$w_{2n}^{-1}$ lies in $S_1(w_{2n}^{-1}w_n, d(1, w_{2n}^{-1}w_n) - K)$
tends to one exponentially quickly. The pair $(w_n, w_{2n}^{-1} w_n)$
is independent, and distributed as $\mun{n} \cross \rmun{n}$.  We may
then use the fact that the measure for shadows of points decays
exponentially in $r$ to show that the $\mun{n} \cross \rmun{n}$
measure of a shadow of the diagonal in $\Gbar \cross \Gbar$ also
decays exponentially in $r$. As $r$ grows linearly in $n$, this shows
that the probability $w_n$ has bounded translation length decays
exponentially in $n$.

\subsection{Acknowledgements} \label{section:acknowledgements}

The author would like to thank G. Margulis and D. Thurston for useful
conversations, and A. Lenhzen for pointing out an error in the proof
of Lemma \ref{lemma:nested sets} in an earlier version of this paper.
The author was supported by NSF grant DMS 0706764.  Support for this
project was also provided by PSC-CUNY Award 60019-40 41, jointly
funded by The Professional Staff Congress and The City University of
New York.

\section{Preliminaries} \label{section:preliminaries}

\subsection{Random walks}

We now review some background on random walks on groups, see for
example Woess \cite{woess}.  Let $G$ be the mapping class group of an
orientable surface of finite type, which is not a sphere with three or
fewer punctures, and let $\mu$ be a probability distribution on $G$.
We may use the probability distribution $\mu$ to generate a Markov
chain, or \emph{random walk} on $G$, with transition probabilities
$p(x,y) = \mu(x^{-1}y)$. We shall always assume that we start at time
zero at the identity element of the group.  The \emph{step space} for
the random walk is the product probability space $(G, \mu)^{\Z_+}$,
and we shall write $(s_1, s_2, \dots)$ for an element of the step
space. The $s_i$ are a sequence of independent, identically
$\mu$-distributed random variables, which we shall refer to as the
\emph{increments} of the random walk. The location of the random walk
at time $n$ is given by $w_n = s_1 s_2 \ldots s_n$, and so the
distribution of random walks at time $n$ is given by the $n$-fold
convolution of $\mu$, which we shall write as $\mun{n}$.  The
\emph{path space} for the random walk is the probability space
$(G^{\Z_+},\P)$, where $G^{\Z_+}$ is the set of all infinite sequences
of elements $G$, and the the measure $\P$ is induced by the map $(s_1,
s_2, \ldots) \mapsto (w_1, w_2, \ldots)$.

We shall always require that the group generated by the support of
$\mu$ is \emph{non-elementary}, which means that it contains a pair of
\pA elements with distinct fixed points in $\PMF$. We do not assume
that the probability distribution $\mu$ is symmetric, so the group
generated by the support of $\mu$ may be strictly larger than the
semi-group generated by the support of $\mu$. Throughout this paper we
will need to assume that the probability distribution $\mu$ has finite
support.

In \cite{maher1}, we showed that it followed from results of
Kaimanovich and Masur \cite{km} and Klarreich \cite{klarreich}, that a
sample path converges almost surely to a uniquely ergodic, and hence
minimal, foliation in the Gromov boundary of the relative space. This
gives a measure $\nu$ on $\fmin$, known as \emph{harmonic measure}.
The harmonic measure $\nu$ is $\mu$-stationary, i.e. \[ \nu(X) =
\sum_{g \in G} \mu(g)\nu(g^{-1}X). \]

\begin{theorem} \cites{km, klarreich, maher1}
\label{theorem:converge}
Consider a random walk on the mapping class group of an orientable
surface of finite type, which is not a sphere with three or fewer
punctures, determined by a probability distribution $\mu$ such that
the group generated by the support of $\mu$ is non-elementary. Then
a sample path $\{ w_n \}$ converges to a uniquely ergodic
foliation in the Gromov boundary $\fmin$ of the relative space $\Grel$
almost surely, and the distribution of limit points on the boundary
is given by a unique $\mu$-stationary measure $\nu$ on $\fmin$.
\end{theorem}

It will also be convenient to consider the \emph{reflected} random
walk, which is the random walk generated by the reflected measure
$\rmu$, where $\rmu(g) = \mu(g^{-1})$. We will write $\rnu$ for the
corresponding $\rmu$-stationary harmonic measure on $\fmin$.

\subsection{Coarse geometry} \mylabel{section:coarse}

We briefly recall some useful facts about Gromov hyperbolic or
$\delta$-hyperbolic spaces, and fix some notation. A
$\delta$-hyperbolic space is a geodesic metric space which satisfies a
\emph{$\delta$-slim triangles} condition, i.e. there is a constant
$\delta$ such that for every geodesic triangle, any side is contained
in a $\delta$-neighbourhood of the other two. Let $(G, d)$ be a
$\delta$-hyperbolic space, which need not be proper.  We shall write
$\d G$ for the Gromov boundary of $G$, and let $\Gbar = G \cup \d G$.
Given a subset $X \subset \Gbar$, we shall write $\overline{X}$ for
the closure of $X$ in $\Gbar$.  Given a point $z \in G$, the Gromov
product based at $z$ is defined to be
\[ \gp{z}{x}{y}  = \tfrac{1}{2} ( d(z, x) + d(z, y) - d(x, y) ). \]
We may extend the definition of the Gromov product to points on
  the boundary by
\[ (x \cdot y)_z  = \sup \liminf_{i,j \to \infty} (x_i \cdot y_i)_z \]
where the supremum is taken over all sequences $x_i \to x$ and $y_j
\to y$. This supremum is finite unless $x$ and $y$ are the same point
in $\partial G$.

We will make use of the following properties of the Gromov product,
see for example, Bridson and Haefliger \cite{bh}*{III.H 3.17}.

\begin{properties}[Properties of the Gromov product] \mylabel{prop:gp}

\ \\

\vspace{-0.2in}

\begin{enumerate}

\item The Gromov product $\gp{z}{x}{y}$ is equal to the distance from $z$ to
a geodesic from $x$ to $y$, up to a bounded error of at most
$\delta$.  \mylabel{prop:geodesic}

\item For any three points $x, y , z \in \Gbar$,
\[ \gp{1}{x}{y} \ge \min \{ \gp{1}{x}{z}, \gp{1}{y}{z} \} - 2
\delta. \] \mylabel{prop:gpmin}

\item If $y \in \d G$, then there is a sequence $y_i \to y$ with
$\lim_n \gp{1}{x}{y_i} = \gp{1}{x}{y}$. \mylabel{prop:seqexist}

\item For any $x \in \Gbar$, and for any sequence $y_i \to y \in \d G$,
\[ \gp{1}{x}{y} - 2 \delta \le \liminf_i \gp{1}{x}{y_i} \le
\gp{1}{x}{y}. \] \mylabel{prop:seq}
\end{enumerate}

\end{properties}

We will also use the following stability property of quasi-geodesics
in a $\delta$-hyperbolic space. Let $I$ be a connected subset of
$\R$. A \emph{quasi-geodesic} is a map $\gamma \colon I \to G$ which
coarsely preserves distance, i.e. there are constants $K$ and $c$ such
that
\[ \frac{1}{K} \norm{s - t} - c \le d(\g(s), \g(t)) \le K \norm{s -
  t} + c. \]
For every $K$ and $c$ there is a constant $L$, which depends only on
$K, c$ and $\delta$, such that a finite $(K, c)$-quasigeodesic is
Hausdorff distance at most $L$ from a geodesic connecting its
endpoints, see Bridson and Haefliger \cite{bh}*{III.H Theorem 1.7}.

\mydef\Knpp{K_1}
Finally, we will also use the fact that nearest point projection onto
a geodesic $\gamma$ is coarsely well defined, i.e. there is a constant
$K$, which only depends on $\delta$, such that if $p$ and $q$ are
nearest points on $\gamma$ to $x$, then $d(p,q) \le K$. Furthermore,
if $y$ is a point on a geodesic $\g$, and $x$ is a point with nearest
point projection $p$ on $\g$, then the path consisting of a geodesic
from $x$ to $p$, and then from $p$ to $y$ is a bounded Hausdorff
distance $\Knpp$ from a geodesic from $x$ to $y$, where $\Knpp$ only
depends on $\delta$, see for example \cite{maher2}*{Proposition 3.1}.

\subsection{Shadows} \mylabel{section:shadows}

Given a point $x \in \Gbar$ and a real number $r$, we define the
\emph{shadow} of $x$ based at $1$, written as $S_1(x, r)$, to be
\[ S_1(x, r) = \{y \in \Gbar \mid \gp{1}{x}{y} \ge r \}. \]
If $x \in G$, and $r \ge d(1, u) + 2 \delta$, then $S_1(x, r)$ is
empty. If $r \le 0$, then $S_1(x, r)$ consists of all of $\Gbar$.

We warn the reader again that this definition of a shadow differs
slightly from that of other authors, for example Blach\`ere,
Ha\"issinsky, and Mathieu \cite{bhm}, who define their shadows
$\mho(1, r)$ to be $S_1(x, d(1, x) - r) \cap \partial G$, in our
notation.  We also remark that it is possible to use the Gromov
product to define a metric on the Gromov boundary, where roughly
speaking the distance between two boundary points is $e^{-\e d}$,
where $d$ is the Gromov product of the two points based at $1$. In
this case, the intersection of a shadow with the boundary is a small
metric neighbourhood of the boundary point. However, we wish our
neighbourhoods to include points in $G$, for which the boundary metric
is not defined, so we find our definition of shadows more convenient.

We may extend the definition of shadows from points to arbitrary
subsets of $\Gbar$.  Given a subset $U \subset \Gbar$, we define the
\emph{shadow} of $U$ based at $1$, written $S_1(U, r)$, to be the
union of the shadows of all points of $U$, i.e.
\[ S_1(U, r) = \bigcup_{x \in U} S_1(x, r). \]
Note that if $U$ contains points in $G$, then in general $U \not
\subset S_r(U)$. However, it is not hard to show that $\bigcap S_1(U,
r) = \overline{U} \cap \d G$, though we will not use this fact
directly.

There is a lower bound on the Gromov product of any two points in the
shadow of a single point, which we now state as a proposition. This is
a direct consequence of Property \ref{prop:gp}.\ref{prop:gpmin} above.

\begin{proposition} \mylabel{prop:product bound} 
For any $y, z \in S_1(x, r)$, the Gromov product $\gp{1}{y}{z}
\ge r - 2 \delta$.
\end{proposition}

Shadows are closed subsets of $\Gbar$, and we now show that a shadow
of a point is the closure of its intersection with $G$.

\begin{proposition} \mylabel{prop:closure} 
$S_1(x, r) = \overline{S_1(x, r) \cap G}$.
\end{proposition}

\begin{proof}
Suppose $y_i$ is a sequence of points in $S_1(x, r) \cap G$, and $y_i
\to y \in \d G$. Then, by the definition of a shadow,
$\gp{1}{x}{y_i} \ge r$ for all $i$. This implies that $\liminf
\gp{1}{x}{y_i} \ge r$, and so $\sup \liminf \gp{1}{x}{y_i} \ge
r$. Therefore, by the definition of the Gromov product for points in
the boundary, $\gp{1}{x}{y} \ge r$, and so $y \in S_1(x, r)$.

Conversely, if $y \in S_1(x, r) \cap \d G$, then by Property
\ref{prop:gp}.\ref{prop:seqexist}, there is a sequence $y_i \to y$ with
$\lim_n \gp{1}{x}{y_i} = \gp{1}{x}{y}$, and as the Gromov product
takes values in $\Z$, we may pass to a subsequence such that
$\gp{1}{x}{y_i} = \gp{1}{x}{y}$, and so this gives a sequence $y_i$
contained in $S_1(x, r)$, which converges to $y$.
\end{proof}

We shall write $\eta_D(T)$ for all points which lie in a metric
$D$-neighbourhood of $T \cap G$, i.e.
\[\eta_D(T) = \{g \in G \mid d(g, t) \le D \text{ for some } t \in T
\cap G \}. \]
We now show that all points in a metric $D$-neighbourhood of a shadow
of $T$, are contained in a slightly larger shadow of $T$.

\begin{proposition} \mylabel{prop:metric nest}
For any $D \ge 0$, 
\[ \eta_D(S_1(T, r)) \subset S_1(T, r - D). \]
\end{proposition}

\begin{proof}
If $g \in \eta_D(S_r(T))$, then there is a point $h \in S_r(T)$ with
$d(g, h) \le D$, and a point $t \in T$ with $\gp{1}{h}{t} \ge r$. By
the definition of the Gromov product, $\gp{1}{g}{t} \ge \gp{1}{h}{t} -
D$ which in turn is at least $r - D$, so $g \in S_{r-D}(T)$, as required.
\end{proof}

We will use the following properties of shadows of points, which
follow from elementary arguments, see Calegari and Maher \cite{cm} for
detailed proofs. We state the results using the current notation of
this paper.

\pic{}{80}{Sufficiently nested shadows are metrically nested}

\mydef\Knest{K_2}
\begin{lemma}[Nested shadows are metrically nested]
  \mylabel{lemma:nested shadows} \cite{cm}*{Lemma 4.5} 
There is a constant $\Knest$, which only depends on $\delta$, such that for
all positive constants $A$ and $r$, and any $x, z \in G$ with $d(x, z)
\ge A + r + 2\Knest$, the shadow $S_z(x, r)$ is disjoint from the
complement of the shadow $S_z(x, r - A - \Knest)$. Furthermore for any pair
of points $a, b \in G$ such that $a \in S_z(x, r)$ and $b \in G
\setminus S_z(x, r - A - \Knest)$, the distance between $a$ and $b$ is at
least $A$.
\end{lemma}

\pic{}{80}{Changing basepoint for a shadow}

\mydef\Kbasea{K_3}
\mydef\Kbaseb{K_4}
\begin{lemma}[Change of basepoint for shadows]
  \mylabel{lemma:change basepoint} \cite{cm}*{Lemma 4.7} 
There are constants $\Kbasea$ and $\Kbaseb$, which only depend on
$\delta$, such that for any $r$, and any three points $x, y, z \in G$ with
$\gp{z}{x}{y} \le r - \Kbasea$, there is an inclusion of
shadows,
\[ S_z(x, r) \subset S_y(x, s), \]
where $s = d(x, y) - d(x, z) + r - \Kbaseb$. 
\end{lemma}

\mydef\Kcomp{K_5}
\begin{lemma}[The complement of a shadow is approximately a shadow]
  \mylabel{lemma:shadow complement} \cite{cm}*{Lemma 4.6} 
There is a constant $\Kcomp$, which only depends on $\delta$, such
that for all constants $r \ge \Kcomp$, and all $x, z \in G$ with
$d(x, z) \ge r + 2 \Kcomp$,
\[ S_x(z, d(x,z) - r + \Kcomp ) \subset G \setminus S_z(x, r) \subset
S_x(z, d(x,z) - r - \Kcomp).  \]
\end{lemma}

We may further extend the definition of a shadow to subsets of $\Gbar
\cross \Gbar$. Let $U \subset \Gbar \cross \Gbar$, and define the
\emph{shadow} $S_1(U, r)$ to be
\[ S_1(U, r) = \bigcup_{(x, y) \in U} S_1(x, r) \cross S_1(y, r). \]
We shall continue to write $S_1(U, r)$ for the shadow in this
case. Hopefully this will not cause confusion, as it should be clear
from context whether $T$ is a subset of $\Gbar$ or $\Gbar \cross
\Gbar$.

Finally, we remark that the lower bound for the Gromov product in a
shadow, \pref{product bound}, immediately implies that the $r$-shadow
of an $s$-shadow is contained in the shadow $S_1(T, \min
\{r,s\}-2\delta)$.

\begin{proposition} \mylabel{prop:nested neighbourhoods}
Let $T$ be a subset of either $\Gbar$ or $\Gbar \cross \Gbar$. Then
\[ S_1(S_1(T, s), r) \subset S_1(T, \min \{r,s\} - 2 \delta), \]
for all $r$ and $s$.
\end{proposition}

\subsection{Exponential decay for shadows}

In this section we show the following upper bounds for measures of
shadows.

\mydef\Kexp{K_6}
\mydef\Kmun{K_7}
\begin{lemma} \mylabel{lemma:exponential shadow}
Let $\mu$ be a finitely supported probability distribution on $G$
whose support generates a non-elementary subgroup, and let $\nu$ be
the corresponding harmonic measure. Then  there are constants $\Kexp$, $\Kmun$
and $c < 1$, such that for any $x$ with $\nhat{x}
\geqslant \Kexp$ and for any $r$,
\[\nu(S_1(x,r)) \leqslant c^{r}, \] 
and
\[ \mun{n}(S_1(x,r)) \le \Kmun c^{r}. \]
The constants $\Kexp$, $\Kmun$ and $c$ depend on $\mu$ and $\delta$,
but not on $r$ or $x$, as long as $d(1, x) \ge \Kexp$.
\end{lemma}

Here we write $\Kmun c^r$ instead of $O(c^r)$, as it will be
convenient to know explicitly the dependence of the implicit constants
in $O(c^r)$. This result also applies to the reflected random walk
generated by the probability distribution $\rmu(g) = \mu(g^{-1})$, and
we may choose the constants to be the same for both random walks.

The proof of this result is essentially the same as the
proof of exponential decay of measures of halfspaces from
\cite{maher2}. Shadows are slightly more general sets than halfspaces,
so the shadow result is not an immediate consequence of the halfspace
result, although the halfspace result does follow from the version for
shadows. Although the shadow version could be deduced using the
halfspace version, this still requires extra work, so we choose to
give an argument here purely in terms of shadows.  We start by giving
some conditions on a family of nested subsets $X_0 \supset X_1 \supset
\cdots$ of $\Gbar$, which guarantee that their measures decay
exponentially in the number of nested sets. We then show that a shadow
$S_1(x, r)$ is contained in a nested family of shadows satisfying the
conditions, and furthermore, the number of sets in the nested family
is linear in $r$.

If $A$ and $B$ are subsets of $\Gbar$, then we define $\dhat{A, B}$,
the distance between $A$ and $B$, to be the smallest distance between
any pair of points in $A \cap G$ and $B \cap G$. If either of these
sets is empty, the distance is undefined.

\begin{lemma} \mylabel{lemma:nested sets}
Let $\mu$ be a probability distribution of finite support of diameter $D$.
Let $X_0 \supset X_1 \supset X_2 \supset \cdots$ be a sequence of
nested closed subsets of $\Gbar$ with the following properties:
\begin{align} 
& 1 \not  \in X_0 \label{eq:p1} \\
& X \setminus X_i \cap X_{i+1} = \varnothing \label{eq:p2} \\
& \dhat{X \setminus X_i, X_{i+1}} \ge D \label{eq:p3} \\
\intertext{%
Furthermore, suppose there is a constant $0 < \e < \half$ such that for
any $x \in X_i \setminus X_{i+1}$, 
}
& \nu_x(X_{i+2}) \le \e, \label{eq:p4} \\
& \nu_x(X \setminus X_{i-1}) \le \e, \label{eq:p5}
\end{align}
then there are constants $c < 1$ and $K$, which only depend on $\e$,
such that $\nu(X_i) \le c^i$ and $\mun{n}(X_i) \le K c^i$.
\end{lemma}

\begin{proof}
By properties \eqref{eq:p1}, \eqref{eq:p2} and Proposition 2.4, any
sequence of points which converges into the limit set of $X_{i+2}$
must contain points in $X_{i+1}$. As the diameter of the support of
$\mu$ is $D$, property \eqref{eq:p3} implies that any sample path
which converges into $X_{i+2}$ must contain at least one point in $X_i
\setminus X_{i+1}$.  Therefore, in order to find an upper bound for
the probability a sample converges into $X_{i+2}$, we can condition on
the location at which the sample path first hits $X_i \setminus
X_{i+1}$. Let $F$ be the (improper) distribution of first hitting
times in $X_i$, i.e. $F(x)$ is equal to the probability that a sample
path first hits $x \in X_i $. This is an improper distribution in
general as $F(X_i) = \sum_{x \in X_i} F(x)$ may be strictly less than
one, as there may be sample paths which never hit $X_i$. As $F$ is
supported on $X_i \setminus X_{i+1}$,
\[ \nu(X_{i+2}) = \sum_{x \in X_i \setminus X_{i+1}}
F(x)\nu_x(X_{i+2}). \]
For all $x \in X_i \setminus X_{i+1}$, there is an upper bound
$\nu_x(X_{i+2}) \le \e$, by property \eqref{eq:p4}, so 
\begin{equation} \label{upper}
\nu(X_i) \le \e F(X_i).
\end{equation}

Not all sample paths which converge to $X_{i-1}$ need to hit $X_i$,
but those that hit $X_i$ and then converge to $X_{i-1}$, give a lower
bound on $\nu(X_{i-1})$, i.e.
\[ \nu(X_{i-1}) \ge \sum_{x \in X_i \setminus X_{i+1}}
F(x)\nu_x(X_{i-1}). \]
By property \eqref{eq:p5}, $\nu_x(X_{i-1}) \ge 1 - \e$, so
\begin{equation} \label{lower}
\nu(X_{i-1}) \ge (1-\e)F(X_i)
\end{equation}
Therefore, combining \eqref{upper} and \eqref{lower}, gives
\[ \frac{ \nu(X_{i+2}) }{ \nu(X_{i-1}) } \le \frac{\e}{1 - \e} < 1,\]
as $\e < \half$. Therefore $\nu(X_i) \le c^i$, where we may choose $c
= \sqrt[3]{\e/(1-\e)}$.

The measure $\nu$ is $\mu$-stationary, and so $\mun{n}$-stationary for
all $n$, i.e.
\[ \nu(X_i) = \sum_{g \in G} \mun{n}(g) \nu_x(X_i). \]
As all terms in the sum are positive, we may discard some
of the terms and the sum will still be bounded above by the upper
bound for $\nu(X_i)$, i.e.
\[ c^i \ge \sum_{g \in X_{i+1} \setminus X_{i+2}} \mun{n}(g)\nu_x(X_i). \]
The measure $\nu_x(X_i)$ is at least $1 - \e$ by \eqref{eq:p5}, which implies
\[ c^i \ge \sum_{g \in X_{i+1} \setminus X_{i+2}} \mun{n}(g)(1 - \e), \]
and we may rewrite this as
\[ c^i \ge (1 - \e) \mun{n}(X_{i+1} \setminus X_{i+2} ). \]
As $X_{i+1} = X_{i+1} \setminus X_{i+2} \cup X_{i+2} \setminus X_{i+3}
\cup \cdots$ this implies

\[ \mun{n}(X_i) \le \frac{1}{1 - \e} \frac{1}{1-c} c^i, \] 
so $\mun{n}(X_i) \le K c^i$, where $1/K = (1 - \e)(1 - c)$. The
constant $K$ only depends on $\e$, as $c$ only depends on $\e$.
\end{proof}

We wish to apply this lemma to shadows of points.  We start by showing
that as the harmonic measure $\nu$ is non-atomic, the harmonic measure
of the shadows of points $S_1(x, r)$ tends to zero as $r$ tends to
infinity, uniformly in $x$.

\mydef\Ke{K_8}
\begin{proposition} \mylabel{prop:epsilon}
For any $\e > 0$ there is a constant $\Ke$, which depends on $\e$ and
$\mu$, such that if $r \ge \Ke$ then $\nu(S_1(x, r)) \le
\e$.
\end{proposition}

\begin{proof}
Suppose not, then there is an $\e > 0$, and a sequence of shadows
$S_1( x_i, r_i )$, with $r_i \to \infty$ such that
$\nu(S_1(x_i, r_i)) \ge \e$.  Let $U_n = \bigcup_{i \ge n} S_1(x_i,
r_i)$, and let $U = \bigcap U_n$, so $U$ consists of all points which
lie in infinitely many $r$-shadows. The sets $U_n$ are decreasing,
i.e. $U_n \supset U_{n+1}$, and $\nu(U_n) \ge \e$ for all $n$, so
$\nu(U) \ge \e$, and so in particular $U$ is non-empty.

Given $\l \in U$, pass to a subsequence, which by abuse of notation we
shall still refer to as $S_1(x_i, r_i)$, such that $\l \in S_1(x_i,
r_i)$ for all $i$. Let $y_i$ be any sequence of points with $y_i \in
S_1(x_i, r_i)$. By \pref{product bound}, $\gp{1}{y_i}{\l} \ge r_i - 2
\delta$ which tends to infinity as $i \to \infty$, which implies that
$y_i \to \l$. But this implies $U = \{ \l \}$, which must have measure
zero, as the measure $\nu$ is non-atomic, which contradicts the fact
that $\nu(U) \ge \e > 0$.
\end{proof}

It will be convenient to choose $\e < \half$, so from now on we will
fix a value of $\Ke$ which ensures that \pref{epsilon} holds for some
$\e$ with $\e < \half$. We now complete the proof of Lemma
\ref{lemma:exponential shadow} by showing that a shadow $S_1(x, r)$
has a nested family of sets $X_n$ satisfying Lemma \ref{lemma:nested
  sets}, where the number of sets is linear in $r$. The constant $L$
in Lemma \ref{lemma:nesting} depends only on $\mu$ and $\delta$, as
does the choice of constant $\e$ from \pref{epsilon}, so the constants
arising from Lemma \ref{lemma:nested sets} will depend only on $\mu$
and $\delta$.

\pic{}{100}{Nested shadows.}

\begin{lemma} \mylabel{lemma:nesting}
For any constant $D$, there is constant $L$, which depends on $\mu$
and $\delta$, with the following properties. For any shadow $S_1(x,
r)$, with $\dhat{1, x} > 2L$, let $N$ be the largest integer such that
$N \le r/L-2$.  Then the sets $X_n = S_1(x, L(n+1))$, for $0 \le n \le
N$, form a sequence of nested sets, which contain $S_1(x, r)$, and
which satisfy properties (1--5) from Lemma \ref{lemma:nested sets}
above.
\end{lemma}

\begin{proof}
Let $L = D + 2\Knest + \Kbasea + \Kcomp + \Ke + 2 \delta$, where $D$ is
the diameter of the support of $\mu$, and the constants $K_i$ are the
constants from Lemmas \ref{lemma:nested shadows}, \ref{lemma:change
  basepoint}, \ref{lemma:shadow complement} and \pref{epsilon}
respectively. We may assume that $L > 0$.  The sets $X_n = S_1(x, Ln)$
are nested, i.e. $X_0 \supset X_1 \supset \cdots$, by the definition
of shadows, and $S_1(x, r) \subset X_n$ for $n \le N \le
r/L - 2$. We now check properties (1--5) from Lemma \ref{lemma:nested
  sets}.

\eqref{eq:p1} The Gromov product $\gp{1}{1}{x} = 0$. For all $y \in X_0$ the
Gromov product $\gp{1}{x}{y} \ge L > 0$, so $1 \not \in X_0$.

\eqref{eq:p2} By Property \ref{prop:gp}.\ref{prop:seq} of the Gromov
product, for any sequence $y_i \to y \in \d G$, $\liminf_i
\gp{1}{x}{y_i} \ge \gp{1}{x}{y} - 2 \delta$.  Therefore, if $y \in
X_{n+1} = S_1(x, L(n+2))$, then for any sequence $y_i \to y$, all but
finitely many points lie in $X_n = S_1(x, L(n+1))$, as $L > 2
\delta$. Therefore $X_{n+1} \cap X \setminus X_n = \varnothing$, as
required.

\eqref{eq:p3} Two shadows which are sufficiently nested in
terms of their shadow parameters, are also metrically nested in terms
of the distance in $G$, by Lemma \ref{lemma:nested shadows}. We shall
apply Lemma \ref{lemma:nested shadows}, choosing the constant $A$ to
be $D$ and the constant $r$ to be $nL$. Recall that $L \ge D +
2\Knest$, where $D$ is the diameter of the support of $\mu$, and
$\Knest$ is the constant from Lemma \ref{lemma:nested shadows}. This
implies that $\dhat{1, x} \ge D + L(n+1) + 2 \Knest$ for all $0 \le n
\le N-1$, by our choice of $N$. Therefore Lemma \ref{lemma:nested
  shadows} implies that $\dhat{ S_1(x, L(n+1)), G \setminus S_1(x,
  Ln)} \ge D$, so $ d(X_{n+1}, G \setminus X_n) \ge D$, as required.

\eqref{eq:p4} Suppose that $y \not \in X_{n+1}$. We wish to show that
that $X_{n+2}$ is contained in a shadow with basepoint $y$,
with a lower bound on the size of its $r$-parameter. This in turn will
give an upper bound on the harmonic measure of the shadow. We
may change the basepoint for the shadows using Lemma
\ref{lemma:change basepoint}, so as long as $\gp{1}{x}{y} \le r -
\Kbasea$, Lemma \ref{lemma:change basepoint} implies that the
shadow $S_1(x, r)$ is contained in $S_y(x, s)$, where
\[ s = r + d(x, y) - d(1, x) - \Kbaseb. \]
As $y \not \in X_{n+1}$, this implies that $ \gp{1}{x}{y} < L(n+1)
$. Therefore choosing $r = L(n+2)$ implies that $\gp{1}{x}{y} < r -
L$, and as we have chosen $L > \Kbasea$, the conditions of Lemma
\ref{lemma:change basepoint} are satisfied.

Therefore $\nu_y(S_y(x, s))$ is an upper bound for
$\nu_y(X_{n+2})$. The harmonic measure $\nu_y(S_y(x, s))$ is equal to
$\nu(S_1(y^{-1}x, s))$, and this is at most $\e < \half$ as long as $s
\ge \Ke$, by \pref{epsilon}. We now verify this last inequality.  By
the definition of the Gromov product,
\[ d(x, y) - d(1, x) = d(1, y) - \gp{1}{x}{y}. \]
As $d(1, y) \ge 0$, and $\gp{1}{x}{y} <
L(n+1)$ this implies that $s \ge L - \Kbasea$. As we have chosen $L >
\Ke + \Kbasea$, this implies that $s \ge \Ke$, as required.

\eqref{eq:p5} Suppose that $y \in X_{n}$. We wish to show that $G
\setminus X_{n-1}$ is contained in a shadow with basepoint $y$, with a
lower bound on the size of its $r$-parameter, which gives an upper
bound on the harmonic measure of the shadow. We have chosen $L$ such
that $L(n-1) \ge \Kcomp$, and $d(1, x) \ge L(n-1) + 2 \Kcomp $, so by
Lemma \ref{lemma:shadow complement},
\[ G \setminus X_{n-1} = G \setminus S_1(x, L(n-1) ) \subset S_x(1, r
), \]
where $r = d(1, x) - L(n-1) - \Kcomp$.  The argument is now
essentially the same as in case (4), except with $1$ and $x$
interchanged. Let $y \in X_n$, so $\gp{1}{x}{y} \ge Ln$.
As $L \ge \Kbasea + \Kcomp$, we may apply Lemma \ref{lemma:change
  basepoint}, which implies that $ S_x(1, r ) \subset S_y(1, s )$,
where 
\[ s = d(1, y)- L(n-1) - \Kcomp - \Kbaseb. \]
We now wish to use \pref{epsilon} to find an upper bound for
$\nu_y(S_y(1, s))$ which is equal to $\nu(S_1(y^{-1}, s))$. By thin
triangles and the definition of the Gromov product, $d(1, y) \ge
\gp{1}{x}{y} - 2 \delta$, so
\[ d(1, y) \ge L(n-1) + L - 2 \delta, \]
which we may rewrite as
\[ d(1, y) - L(n-1) - \Kbaseb - \Kcomp \ge L -\Kbaseb - \Kcomp - 2
\delta, \]
where the left hand side is equal to $s$. As we have chosen $L \ge
\Kbaseb + \Kcomp + \Ke + 2 \delta $, this shows that $s \ge
\Ke$. Therefore \pref{epsilon} implies that
$\nu(S_1(y^{-1}, s)) \le \e$, so $\nu_y(X_{n-1}) \le \e < \half$, as
required.
\end{proof}

\section{Linear progress} \label{section:linear progress}

In this section we prove Theorem \ref{theorem:exponential linear
  progress}, i.e. we show that sample paths make linear progress at
some rate $L$, and furthermore, the proportion of sample paths at time
$n$ which are distance at most $Ln$ from $1$ decays exponentially in
$n$. As $d(1, g)$ is equal to $d(1, g^{-1})$ the reflected random walk
also makes linear progress at the same rate $L$, and with the same
exponential decay constant for the proportion of sample paths distance
less than $Ln$ from the origin $1$.

A random walk of length $nk$, determined by a probability distribution
$\mu$, may be thought of as a random walk of length $n$, determined by
the probability distribution $\mu_k$. We shall write $w^k_n$ for
$w_{kn}$, and we shall call this the \emph{$k$-iterated random
  walk}. The steps of the $k$-iterated random walk are $s^k_i =
s_{(i-1)k+1} \ldots s_{ik}$, and so for each $i$, the segment of the
random walk from $w^k_i$ to $w^k_{i+1}$ is independently and
identically distributed according to the probability distribution
$\mu_k$, the $k$-fold convolution of $\mu$. However, the distance from
$1$ to $w^k_{i+1}$ is at most $\dhat{1, w^k_i} + \dhat{w^k_i,
  w^k_{i+1}}$, but may be smaller, as the random walk may have
``backtracking,'' i.e. the geodesic from $w^k_i$ to $w^k_{i+1}$ may
fellow travel with a terminal segment of the geodesic from $1$ to
$w^k_i$. This is illustrated schematically in Figure
\ref{pic:backtracks}, for the first few steps of the $k$-iterated
random walk.

Set $X^k_i$ to be the random variable corresponding to the change in
distance from the basepoint $1$ from time $i-1$ to time $i$ of the
$k$-iterated random walk, i.e.
\[ X^k_i = \dhat{1, w^k_{i}} - \dhat{1, w^k_{i-1} }, \]
which may be negative. The sum of the first $n$ random variables
$X^k_i$ is equal to the distance travelled at the $n$-th step of the
$k$-iterated walk, i.e.
\[ \sum_{i = 1}^{n} X^k_i = \dhat{1, w^k_n}. \]
We may write $X^k_i = Y^k_i - Z^k_i$, where $Y^k_i$ is the distance
the $k$-iterated random walk travels between steps $i-1$ and $i$, i.e.
\[ Y^k_i = \dhat{w_{(i-1)k}, w_{ik}}, \]
and $Z^k_i = Y^k_i - X^k_i$. The $Y^k_i$ form an independent
collection of random variables, but the $Z^k_i$ do not. By the
definition of the Gromov product,
\[ Z^k_i = 2\gp{w^k_{i-1}}{1}{w^k_i}, \]
and we may think of $Z^k_i$ as the amount of backtracking the iterated
random walk $w^k_n$ does from step $i-1$ to step $i$. In particular,
the $Z^k_i$ are non-negative. In order to find lower bound estimates
for the sums of the $X_i$, it suffices to find lower bound estimates
for the sums of the $Y_i$, and upper bound estimates for the sums of
the $Z_i$, and we now show how to do this, using standard results from
the theory of concentration of measures.

The distances $Y^k_i = \dhat{w^k_{i-1}, w^k_i}$ form a sequence of
independent, identically distributed random variables, so estimates on
the behaviour of the sums of these random variables are well
known. Let $Y^k$ be the expected value of $Y^k_i$, which depends on
$k$, but not on $i$. As the trajectories of the random walk converge
to the boundary almost surely, $Y^k \to \infty$ as $k \to \infty$.  We
will use the following Bernstein or Chernoff-Hoeffding estimate, see
for example Dubhashi and Panconesi \cite{dp}*{Theorem 1.1} which says
that the probability that the sum of $n$ copies of $Y^k_i$ deviates
from the expected mean $nY^k$ by at least $\e n$ decays exponentially
in $n$.

\begin{theorem} \mylabel{theorem:bernstein}
Let $Y_i$ be a sequence of bounded independent identically distributed random
variables with mean $Y$.  Then for any $\e > 0$ there is a constant $c
< 1$ such that
  \[ \P \left( \norm{ \sum_{i=1}^n (Y_i - Y) } \ge \e n \right) \le c^{n}. \]
\end{theorem}

We now show a similar bound for the sums of the $Z^k_i$. We start by
showing that the distribution functions of the $Z^k_i$ are bounded
above by the same exponential function, for all $k$ and
$i$. Furthermore, the upper bound for $Z^k_i$ holds independently of
the values of $Z^k_j$ for $j < i$. As $Z^k_i$ is a function of $w^k_j$
for $j \le i$, it suffices to show that the upper bound is independent
of the values of $w^k_j$ for $j < i$.

\begin{proposition} \label{prop:gp estimate} 
There are constants $K$ and $c < 1$, which do not depend on $k$ or
$i$, such that
\[ \P(Z^k_i \ge r \mid w^k_1, \ldots, w^k_{i-1}) \le K c^r. \]
\end{proposition}

\begin{proof}
By the definition of $Z^k_i$, if $Z^k_i \ge r$ then $\gp{ w^k_{i-1} }{
  1 }{ w^k_i } \ge \half r$. By the definition of shadows, this
condition is equivalent to the condition $w^k_i \in S_{w^k_{i-1}}(1,
\tfrac{1}{2}r )$, therefore
\[ \P(Z^k_i \ge r \mid w^k_1, \ldots, w^k_{i-1} ) = \P \left( w^k_i
  \in S_{w^k_{i-1}}( 1 , \tfrac{1}{2}r ) \mid w^k_1, \ldots, w^k_{i-1} \right)
. \]
We may apply the isometry $(w^k_{i-1})^{-1}$, and use the fact that
$(w^k_{i-1})^{-1}w^k_i = s^k_i$, to obtain,
\[ \P(Z^k_i \ge r \mid w^k_1, \ldots, w^k_{i-1} ) = \P \left( s^k_i
  \in S_1( (w^k_{i-1})^{-1}, \tfrac{1}{2} r ) \mid w^k_1, \ldots,
  w^k_{i-1} \right) . \]
As $s^k_i$ is distributed as $\mun{k}$, and is independent of the
$w^k_j$ for $j < i$, this implies,
\[ \P(Z^k_i \ge r \mid w^k_1, \ldots, w^k_{i-1} ) = \mun{k} ( S_1(
(w^k_{i-1})^{-1}, \tfrac{1}{2} r ). \]
Now using Lemma \ref{lemma:exponential shadow}, there are constants
$\Kmun$ and $c < 1$ such that the bound $\mun{k}(S_1(g, \tfrac{1}{2}
r)) \le \Kmun c^{r/2}$, is independent of $g$ and $k$, so this implies
\[ \P(Z^k_i \ge r \mid w^k_1, \ldots, w^k_{i-1} ) \le \Ke c^{r/2}, \]
as required.
\end{proof}

In particular, this gives an upper bound for the expected value of
$Z^k_i$ which is independent of $k$. Therefore, by choosing $k$ to be
large, we can make the expected value of $Y^k_i$ much larger than the
expected value of $Z^k_i$.

We now show that there is a constant $L >0$, which is independent of
$k$, such that the probability that the sum $Z^k_1 + \cdots + Z^k_n$
is larger than $Ln$ decays exponentially in $n$.

\begin{lemma} \mylabel{lemma:Z estimate}
Let $w^k_n$ be the $k$-iterated random walk of length $n$, generated by a
finitely supported probability distribution $\mu$, whose support generates a
non-elementary subgroup of the mapping class group, and let $Z^k_i = 2
\gp{w^k_{i-1}}{1}{w^k_i}$. Then there are constants $L, K$ and $c < 1$,
which depend on $\mu$ but are independent of $k$, such that
\[ \P(Z^k_1 + \cdots + Z^k_n \ge Ln ) \le K c^n, \]
for all $n$.
\end{lemma}

\begin{proof}
We have shown that the probability that $Z^k_i \ge r$ decays
exponentially in $r$, with exponential decay constants which do not
depend on either $k$ or $i$, or the values of any other $Z^k_j$ for $j
< i$. As $Z^k_i$ is also independent of $Z^k_j$ for $j > i$, this
implies that the exponential bounds for $Z^k_i$ hold independently of
the vales of $Z^k_j$ for all $j \not = i$. Therefore, the probability
distribution of the sum $Z^k_1 + \cdots + Z^k_n$ will be bounded above
by a multiple $K^n$ of the $n$-fold convolution of the exponential
distribution function with itself. We will use the following
Chernoff-Hoeffding bound for sums of exponential random variables. The
version stated below is an exercise from Dubhashi and Panconesi
\cite{dp}, but we provide a proof in Appendix \ref{appendix} for
completeness.

\begin{proposition} \cite{dp}*{Problem 1.10}.
Let $A_i$ be independent identically distributed exponential random
variables, with expected value $A$. Then for any $t \ge 0$,
\[ \P( A_1  + \cdots + A_n  \ge ( 1+ t ) nA ) \le \left( \frac{1 +
    t}{e^t} \right)^n. \]
\end{proposition}

The upper bound for the sum of the $Z^k_i$ will therefore be $K^n$
times the upper bound above, i.e.
\[ \P( Z^k_1 + \cdots + Z^k_n \ge ( 1 + t ) n Z^k) \le \left( \frac{1
    + t}{e^t} K \right)^n.  \]
The expected value $Z^k$ is bounded above for all $k$, so by choosing
$t$ sufficiently large, we may ensure that the base of the exponent
on the right is strictly less that $1$.  This completes the proof of
Lemma \ref{lemma:Z estimate}.
\end{proof}

We now complete the proof of Theorem \ref{theorem:exponential linear
  progress}. Recall that $d(1, w^k_n) = X^k_1 + \cdots + X^k_n$, and
$X^k_i = Y^k_i - Z^k_i$. So if
\[ \sum_{i=1}^n X^k_i \le n(Y^k - \e - ( 1 + t )  Z^k), \] 
then $Y^k_1 + \cdots + Y^k_n - n Y^k \le - \e n $, or $Z^k_1 +
\cdots + Z^k_n \ge ( 1 + t ) n Z^k$, though of course both
conditions may be satisfied. Furthermore, we may choose $k$
sufficiently large such that $L = Y^k - \e - (1 + t) Z^k$ is
positive. The probability that at least one of the events occurs is at
most the sum of the probability that either occurs, so
\[ \P(d(1, w^k_n) \le Ln) \le c_1^n + K c_2^n, \]
for some constants $c_1 < 1$ from Theorem \ref{theorem:bernstein} and
$c_2 < 1$ from Lemma \ref{lemma:Z estimate}, and this decays
exponentially in $n$, as required. This completes the proof of Theorem
\ref{theorem:exponential linear progress}.

\section{Translation length} \label{section:translation length}

In this section we prove Theorem \ref{theorem:main}. We start by
showing that translation length of $g$ is coarsely equivalent to the
length of the (relative) shortest element in the conjugacy class of
$g$, which we shall denote $[g]$, i.e.
\[ [g] = \inf_{h \in G} d(1, hgh^{-1}). \]

\begin{lemma}
Let $G$ be the mapping class group of a non-sporadic surface.
There is a constant $K$ such that $\norm{ \t(g) - [g]} \le K$.
\end{lemma}

\begin{proof}
Let $g$ be a conjugate of minimal relative length $[g]$.  By the
definition of translation length, $\t(g) \le \dhat{1, g} = [g]$. We
now show the bound in the other direction.

There is a constant $M$, which depends only on the surface, such that
every non-\pA element is conjugate to an element of relative length at
most $M$, see for example \cite{maher1}*{Lemma 5.5} so we shall choose
$K > M$ and then we may assume that $g$ is \pA.  Let $\a$ be a
quasi-axis for $g$, i.e. a bi-infinite quasigeodesic such that $\a$
and $g^n \a$ are $2 \delta$-fellow travellers for all $n$. Let $h$ be
a closest point on $\a$ to $1$, then $gh$ is distance at most $\t(g) +
K$ from $\a$. This implies that the distance from $h$ to $gh$ is at
most $\t(g) + K$, so $\dhat{1, h^{-1}gh}$ is at most $\t(g) +
K$. Therefore $[g] \le \t(g) + K$, as required.
\end{proof}

The mapping class group has \emph{relative conjugacy bounds},
\cite{maher1}*{Theorem 3.1}, i.e. there is a constant $K$, which only
depends on the surface, such that if $a$ and $b$ are conjugate, then
$a = vbv^{-1}$ for some element $v$ with
\[ \nhat{v} \le K(\nhat{a} + \nhat{b}). \]
If $g$ is a group element, then we may think of $g$ as a point in the
metric space $(G, d)$. However, we can also represent $g$ by a choice
of geodesic in $G$ from $1$ to $g$. Geodesics need not be unique, but
any two distinct choices of geodesics are Hausdorff distance at most
$2 \delta$ apart. This gives two ways of representing a product $gh$
of two group elements $g$ and $h$. We may choose a single geodesic
from $1$ to $gh$, or alternatively choose a path from $1$ to $gh$
consisting of two geodesic segments, the first consisting of a
geodesic from $1$ to $g$, and the second consisting of a geodesic from
$g$ to $gh$, which is the translate of a geodesic from $1$ to
$h$. Therefore if $g$ is equal to $vsv^{-1}$, we can represent $g$ by
a path composed of three geodesic segments, each consisting of a
translate of $v$, $s$ and $v^{-1}$ respectively, and this is what we
mean when we refer to the \emph{path} $vsv^{-1}$.  The fact that $G$
has relative conjugacy bounds implies that if an element $g$ is
conjugate to a short element $s$, and $v$ is a conjugating word of
shortest possible relative length, then the path $vsv^{-1}$ is
quasigeodesic, where the quasigeodesic constants depend on $\nhat{s}$,
the constant of hyperbolicity $\delta$, and the relative conjugacy
bounds constant $K$.

\begin{lemma} \cite{maher1}*{Lemma 4.2}
Let $G$ be a weakly relatively hyperbolic group with relative
conjugacy bounds. Let $g$ be an element of $G$ which is conjugate to
an element $s$, i.e. $g = v s v^{-1}$, for some $v \in G$. If we
choose $v$ to be a conjugating word of shortest relative length, then
the word $vsv^{-1}$ is quasi-geodesic in the relative metric, with quasi-geodesic
constants which depend only on the relative length of $s$, and the
group constants $\delta$ and $K$. 
\end{lemma}

\pic{}{100}{A quasigeodesic path}

\mydef\Kpath{K_9}

\begin{proposition} \mylabel{prop:kpath}
For any constant $T$ there is a constant $\Kpath$, which only depends
on $T$, the constant of hyperbolicity $\delta$, and the relative conjugacy
bounds constant, such that if $g$ is conjugate to an element $s$ of
relative length at most $T$, then $g = vsv^{-1}$ for some $v$ with the
following properties:

\begin{enumerate}

\item $ d(1, v) \ge \tfrac{1}{2} d(1, g) - \Kpath $

\item $g \in S_1(v, d(1, v) - \Kpath)$

\item $1 \in S_g(gv, d(1, v) - \Kpath)$

\end{enumerate}

\end{proposition}

\begin{proof}
Let $g = vsv^{-1}$, where $d(1, s) \le T$ and $v$ is a conjugating
element of shortest (relative) length. The first inequality follows
from the triangle inequality, which implies that $d(1, v) \ge
\tfrac{1}{2} d(1, g) - T/2$. As the path $vsv^{-1}$ is a
quasigeodesic, there is a constant $L$, which only depends on $T$, the
constant of hyperbolicity $\delta$, and the conjugacy bounds constant,
such that the distance from $v$ to a geodesic from $1$ to $g$ is at
most $L$. This implies that if $p$ is the nearest point projection of
$v$ to a geodesic from $1$ to $g$, then $d(1, p) \ge d(1, v) - L$. As
any geodesic from $v$ to $g$ is contained in a $\Knpp$-neighbourhood of
the nearest point projection path, consisting of a geodesic from $v$
to $p$, and then from $p$ to $g$, where $\Knpp$ only depends on
$\delta$. This implies that the distance from $1$ to any geodesic from
$v$ to $g$ is at least $d(1, v) - L - \Knpp$. Finally, as the Gromov
product $\gp{1}{v}{g}$ is equal to the distance from $1$ to a geodesic
from $v$ to $g$, up to bounded additive error $2 \delta$, this implies
that $\gp{1}{v}{g} \ge d(1, v) - K$, where $K = L + \Knpp + 2
\delta$, which only depends on the constant of hyperbolicity
$\delta$. This means that $g \in S_1(v, d(1, v) - K)$, and the
same argument applied to the points $1, g$ and $vs$ implies that $1
\in S_g(gv, d(1, v) - K)$, as $vs = gv$, for the same constant
$K$. We may therefore choose $\Kpath$ to be the maximum of $K$ and $T/2$.
\end{proof}

\pref{kpath} above shows that the probability that a random walk $w_n$
is conjugate to an element of relative length at most $T$, is bounded
above by the probability that there is an element $v$, with $d(1, v)
\ge \half d(1, w_n) - \Kpath$, such that $w_n \in S_1(v, d(1, v) -
\Kpath)$, and $w_n^{-1} \in S_1(v, d(1, v) - \Kpath)$.  We shall write
$X_n$ for the measure corresponding to the distribution of pairs
$(w_n, w_n^{-1})$ on $\Gbar \cross \Gbar$, i.e.
\[X_n(U) = \P( (w_n, w_n^{-1}) \in U ), \]
for any subset $U \subset \Gbar \cross \Gbar$. As $\P(d(1, w_n)) \le
Ln$ decays exponentially, by Theorem \ref{theorem:exponential linear
  progress}, this gives the following upper bound for that the
probability that $w_n$ is conjugate to an element of length at most
$T$,
\[ P(\tau(w_n) \le T) \le X_n(S_1(\D, \tfrac{1}{2} Ln - \Kpath )) + O(c^n), \]
where $c < 1$ is the constant from Theorem \ref{theorem:exponential
  linear progress}, and where $\D$ is the diagonal in $\Gbar \cross
\Gbar$.

Therefore, in order to complete the proof of Theorem
\ref{theorem:main}, it suffices to show:

\begin{lemma} \mylabel{lemma:Xdiagonal}
Let $L$ be a constant such that $\P(d(1, w_n) \le Ln)$ decays
exponentially in $n$. Then for any $K$, there is a constant $c < 1$,
which depends on $K$ and $\mu$, such that
\[ X_n(S_1(\D, \tfrac{1}{2} Ln - K )) \le O(c^n). \]
\end{lemma}

The rest of this section is devoted to the proof of Lemma
\ref{lemma:Xdiagonal}. In fact, it will be convenient to obtain upper
bounds for $X_{2n}$ rather than $X_n$. This suffices to obtain upper
bounds for $X_n$ for all $n$, as if $D$ is the diameter of the support
of $\mu$, then $X_{2n-1}(S_1(U, r)) \le X_{2n}(S_1(U, r - D))$, by
\pref{metric nest}.

We start by showing that it is very likely that a random walk $w_{2n}$
lies in the shadow $S_1(w_n, \half d(1, w_n))$.

\begin{proposition} 
\mylabel{lemma:walk in halfspace} 
The probability that $w_{2n}$ lies in $S_1(w_n, \half d(1, w_n))$
tends to one exponentially quickly as $n$ tends to infinity, i.e.
\[ \P( w_{2n} \not \in S_1(w_n, \tfrac{1}{2} d(1, w_n)) ) \le O(c^n), \]
for some $c < 1$.
\end{proposition}

\begin{proof}
We shall find an upper bound for the probability that $w_{2n}$ does
not lie in the shadow $S_1(w_n, \tfrac{1}{2}d(1, w_n))$. Conditioning
on $w_n = g$, and using the fact that the complement of the shadow
$S_1(w_n, \tfrac{1}{2}d(1, w_n))$ is contained in $S_{w_n}(1,
\tfrac{1}{2}d(1, w_n) - \Kcomp)$, Lemma \ref{lemma:shadow complement},
gives 
\begin{align*}
\P(w_{2n} \not \in S_1(w_n, \tfrac{1}{2} d(1, w_n)) & \le \sum_{g \in G}
\mun{n}(g) \P(g s_{n+1} \ldots s_{2n} \in S_g(1, \tfrac{1}{2}d(1,
g) - \Kcomp) \mid w_n = g).
\intertext{%
The condition $gs_{n+1} \ldots s_{2n} \in S_g(1, \tfrac{1}{2}d(1, g) -
\Kcomp)$ is the same as $s_{n+1} \ldots s_{2n} \in S_{1}(g^{-1},
\tfrac{1}{2}d(1, g) - \Kcomp)$, and as the $s_{n+1}, \ldots, s_{2n}$ are independent
of $w_{n}$, this implies that
}
\P(w_{2n} \not \in S_1(w_n, \tfrac{1}{2} d(1, w_n)) & \le \sum_{g \in G}
\mun{n}(g) \mun{n}(S_1(g^{-1}, \tfrac{1}{2} d(1,g) - \Kcomp ) ).
\intertext{%
By Theorem \ref{theorem:exponential linear progress}, the probability
that $d(1, w_n) \le Ln$ is at most $O(c_1^n)$, for some $c_1 <
1$, which gives }
\P(w_{2n} \not \in S_1(w_n, \tfrac{1}{2} d(1, w_n)) & \le O(c_1^n) +
\sum_{g \in G \setminus B(1, Ln)} \mun{n}(g) \mun{n}(S_1(g^{-1},
\tfrac{1}{2} d(1,g) - \Kcomp ) ).
\intertext{%
The upper bound for the measure of a shadow, Lemma
\ref{lemma:exponential shadow}, then gives the following upper bound,
}
\P(w_{2n} \not \in S_1(w_n, \tfrac{1}{2} d(1, w_n)) & \le O(c_1^n) +
\sum_{g \in G \setminus B(1, Ln)} \mun{n}(g) O( c_2^{Ln/2 -
  \Kcomp} ),
\intertext{%
for some constant $c_2 < 1$. Therefore
}
\P(w_{2n} \not \in S_1(w_n, \tfrac{1}{2} d(1, w_n)) & \le O(c_1^n) +
O(c_2^{Ln/2}),
\end{align*}
which decays exponentially in $n$, as required.
\end{proof}

Applying this result to the reflected random walk implies that the
probability that $w_{2n}^{-1}$ does not lie in $S_1(w_{2n}^{-1}w_n,
\tfrac{1}{2} d(1, w_{2n}^{-1}w_n))$ also decays exponentially.

We now use this to find an upper bound for $X_{2n}$ in terms of
$\mun{n} \cross \rmun{n}$.

\begin{proposition} \mylabel{lemma:Xmu estimate}
Let $T$ be a subset of $\Gbar \cross \Gbar$. There are constants $L > 0$
and $c < 1$ such that
\[ X_{2n}(S_1(T, r) ) \le \mun{n} \cross \rmun{n}(S_1( T, {\min \{r,
  \tfrac{1}{2} Ln \} - 2 \delta})) + O(c^n). \]
\end{proposition}

\begin{proof}
We have shown that the probability that each of the following four
events occurs tends to one exponentially quickly.
\begin{align*}
& \dhat{1, w_n} \ge Ln \\
& \dhat{1, w_{2n}^{-1}w_n} \ge Ln \\
& w_{2n} \in S_1(w_n, \tfrac{1}{2} d(1, w_n)) \\
& w_{2n}^{-1} \in S_1(w_{2n}^{-1}w_n, \tfrac{1}{2} d(1, w_{2n}^{-1}w_n) )
\end{align*}
Therefore the probability that all four of them occur tends to one
exponentially quickly.

If all four events occur, then $\gp{1}{w_n}{w_{2n}} \ge \half d(1,
w_n) \ge \half Ln$, and similarly,
$\gp{1}{w_{2n}^{-1}w_n}{w_{2n}^{-1}} \ge \half d(1, w_{2n}^{-1}w_n)
\ge \half Ln$. Furthermore, if $(w_{2n}, w_{2n}^{-1})$ lies in
$S_r(T)$, then there is a point $(s, t) \in T$ such that
$\gp{1}{w_{2n}}{s} \ge r$ and $\gp{1}{w_{2n}^{-1}}{t} \ge r$. This
implies that $\gp{1}{w_n}{s} \ge \min \{ r, \half Ln \} - 2 \delta$,
and $\gp{1}{w_{2n}^{-1}w_n}{t} \ge \min \{r, \half Ln \} - 2 \delta$,
and so $(w_n, w_{2n}^{-1}w_n) \in S_1(T, {\min \{r, \half Ln \} - 2
  \delta})$, as required.
\end{proof}

Finally, we now show that the $\mun{n} \cross \rmun{n}$-measure of a
shadow of the diagonal $S_1(\Delta, r)$ decays exponentially in
$r$.

\begin{proposition} \mylabel{lemma:diagonal}
There are constants $c_1 < 1$ and $c_2 < 1$ such that
\[ \mun{n} \cross \rmun{n} (S_1(\Delta, r)) \le O(c_1^r) +
O(c_2^n), \]
for all $n$ and $r$.
\end{proposition}

\begin{proof}
Let $v_n$ and $w_n$ be random walks determined by $\mu$ and $\rmu$
respectively. If $(v_n, w_n) \in S_1(\D, r)$, then there is a point
$x$ such that $\gp{1}{v_n}{x} \ge r$ and $\gp{1}{w_n}{x} \ge
r$. Therefore $\gp{1}{v_n}{w_n} \ge r - 2 \delta$, and so $v_n \in
S_1(w_n, r - 2 \delta)$. By the upper bound for measures of shadows,
Lemma \ref{lemma:exponential shadow}, for any $w_n$ with $d(1, w_n)
\ge \Kexp$, the probability that $v_n \in S_1(w_n, r - 2 \delta)$ is
at most $\Kmun c_1^{r - 2 \delta}$, for some $c_1 < 1$. Furthermore,
by Theorem \ref{theorem:exponential linear progress}, there is a $c_2
< 1$ such that the probability that $d(1, w_n) \le \Kexp$ is at most
$\Kmun c_2^n$, for $n \ge \Kexp / L$. Therefore $\mun{n} \cross
\rmun{n}(S_1(\D, r)) \le O(c_1^r) + O(c_2^n)$, as required.
\end{proof}

Combining Propositions \ref{lemma:Xmu estimate} and \ref{lemma:diagonal}
establishes Lemma \ref{lemma:Xdiagonal}, and so completes the proof of
Theorem \ref{theorem:main}.

\appendix

\section{Chernoff-Hoeffding bounds for exponential random variables}
\mylabel{appendix}

In this section we provide the details for the following
Chernoff-Hoeffding bound for exponential random variables. This proof
is the solution given by Dubhashi and Panconesi to \cite{dp}*{Problem
  1.10}, which appeared in the initial draft version, but not in the
final published version, and we reproduce it here for the sake of
completeness.

\begin{proposition} 
Let $Z_i$ be independent identically distributed exponential random
variables, with expected value $Z$. Then for any $t \ge 0$,
\[ \P( Z_1  + \cdots + Z_n  \ge ( 1+ t ) nZ ) \le \left( \frac{1 +
    t}{e^t} \right)^n. \]
\end{proposition}

\begin{proof}
Let $Z_i$ have probability density function $f(x) = \a e^{-\a x}$, so
the expected value of $Z_i$ is $Z = 1 / \a$, and set $S_n = Z_1 +
\cdots + Z_n$. Consider the moment generating function
\[ \E( e^{\l Z_i}) = \a \int_0^\infty e^{\l x} e^{-\a x} dx  =
\frac{\a}{\a - \l}, \]
for $0 < \l < \a$. Therefore
\[ \E(e^{\l S_n}) = \left( \frac{\a}{\a - \l} \right)^n. \]
It now follows from Markov's inequality that
\[ \P( S_n \ge s ) \le \frac{\E(e^{\l S_n})}{e^{\l s}} = \frac{1}{e^{\l
    s} (1 - \frac{\l}{\a})^n} .\]
The right hand side above is minimized by choosing $\l = \a - n/s$,
which gives
\[ \P( S_n \ge s ) \le \left( \frac{\a s}{n} \right)^n e^{-\a s + n}. \]
Setting $s = (1 + t)n Z$, and using the fact that $Z = 1 / \a$, yields,
\[ \P(S_n \ge (1 + t)nZ ) \le \left( \frac{1 + t}{ e^t } \right)^n, \]
as required.
\end{proof}



\end{document}